\documentclass[10pt]{amsart}

\usepackage{geometry,graphicx,amssymb,amsmath,amsbsy,eucal,amsfonts,mathrsfs,amscd,bm,tcolorbox}

\numberwithin{equation}{section}

\allowdisplaybreaks[3]

\newtheorem{theorem}{Theorem}[section]
\newtheorem{lemma}[theorem]{Lemma}
\newtheorem{assumption}[theorem]{Assumption}
\newtheorem{corollary}[theorem]{Corollary}
\newtheorem{proposition}[theorem]{Proposition}
\theoremstyle{definition}
\newtheorem{definition}[theorem]{Definition}
\theoremstyle{remark}
\newtheorem{remark}[theorem]{Remark}

\newcommand{\R}{\mathbb{R}}

\newcommand{\del}{{\Delta}}
\newcommand{\de}{\delta}

\newcommand{\bl}{\bigl\langle}
\newcommand{\br}{\bigr\rangle}
\newcommand{\bll}{\langle \hspace{-.1cm} \langle}
\newcommand{\brr}{\rangle \hspace{-.1cm} \rangle}

\newcommand{\Th}{\mathcal{T}_h}

\newcommand{\La}{\Lambda}
\newcommand{\la}{\lambda}
\newcommand{\w}{\wedge}

\newcommand{\Q}{\mathsf{P}}

\newcommand{\dd}{\mathsf{d}}

\newcommand\supp{\operatorname{supp}}

\newcommand{\pol}{\EuScript{P}}

\newcommand{\sig}{\sigma}

\newcommand{\C}{\mathcal{C}}

\newcommand{\dr}{\mathcal{R}}

\newcommand{\ZZ}{\mathsf{Z}}
\newcommand{\Zz}{\mathfrak{Z}}

\newcommand{\vol}{\text{{vol}}}
\newcommand{\st}{\operatorname{st}}
\newcommand{\es}{\operatorname{es}}
\newcommand{\bb}{\mathfrak{b}}
\newcommand{\BBb}{\mathsf{b}}

\newcommand{\polb}{\breve{\pol}}

\newcommand{\tr}{\mathsf{t}\mathsf{r}} 
 \newcommand{\sprk}{\pol_r^-\La^k(\Th)}
  \newcommand{\Sp}[2]{\pol_{#1}^-\La^{#2}(\Th)}
 
 \newcommand{\ssp}[3]{\pol_{#1}^-\La^{#2}(#3)}
 \newcommand{\Prk}[2]{P_{#1}^{#2}}

\newcommand{\Zp}{ \mathfrak{Z} \mathring{\pol}_r^{-} \La}
\newcommand{\Zpp}{ \mathfrak{Z}^{\perp} \mathring{\pol}_r^{-} \La}
\newcommand{\ZZp}{ \mathfrak{Z} \polb_r^{-} \La}
\newcommand{\ZZpp}{ \mathfrak{Z}^{\perp} \polb_r^{-} \La}

\newcommand{\UU}{\mathsf{U}}

\title[]{Local $L^2$-bounded commuting projections in FEEC}

\author[]{Douglas Arnold}
\address{Department of Mathematics
University of Minnesota
Minneapolis, MN 55442, USA
}
\email{arnold@umn.edu}

\author[]{Johnny Guzm\'an}
\address{Division of Applied Mathematics
Brown University
Box F
182 George Street
Providence, RI 02912, USA}
\email{johnny\_guzman@brown.edu}

\thanks{The work of the first author was supported by Simons Foundation grant 601937, DNA. The second author was partially supported by NSF grant  DMS-1913083}


\begin{document}

\begin{abstract}
We construct local projections into canonical finite element spaces that appear in the finite element exterior calculus. These projections are bounded in $L^2$ and commute with the exterior derivative.   
\end{abstract}

\maketitle

\section{Introduction}
Bounded commuting projections are a primary instrument in the finite element exterior calculus (FEEC) \cite{arnold2018finite, arnold2006finite, arnold2010finite}.  In particular, the existence of such projections from the Hilbert variant of the de~Rham complex to a finite dimensional
subcomplex is the primary requirement for stable Galerkin approximations for the Hodge Laplacian  \cite[Theorem 3.8]{arnold2010finite}. Error estimates for the Galerkin approximation then follow.  For this, the projections must be bounded on the space of $L^2$ differential forms with exterior derivative in $L^2$. A stronger condition is that the projections are bounded on the larger space $L^2$, which is a primary
requirement to obtain improved error estimates \cite[Theorem 3.11]{arnold2010finite} and also to convergence for eigenvalue problems  \cite[Theorem 3.19]{arnold2010finite}.

The first commuting projections were  developed by Sch\"oberl \cite{schoberlmultilevel} and later by Christiansen and Winther \cite{christiansen2008smoothed}, who treated non-quasiuniform meshes and spaces with essential boundary conditions. In particular, the projection \cite{christiansen2008smoothed} is bounded in $L^2$. However, the projections  \cite{schoberlmultilevel, christiansen2008smoothed}  are not local, meaning that the projection of a form $u$ on a simplex $T$ does not depend solely on $u$ on a patch of elements surrounding $T$. More recently, Falk and Winther \cite{falk2014local, falk2015double} constructed commuting projections that are local and defined for $L^2$ bounded $k$-forms with its exterior derivative also belonging to $L^2$.  More recently in two and three dimensions, Ern et al.  \cite{ern2019equivalence} have developed commuting projections for the last part of the de~Rham complex that are local and bounded in $L^2$. Bounded commuting projections that are local may also be used to obtain error estimates in other norms such as $L^\infty$ \cite{gastaldi1989sharp, demlow2004localized}.

In this paper, inspired by the techniques of Falk and Winther \cite{falk2014local}, we construct local, $L^2$-bounded, commuting projections from the de~Rham complex in any space dimension onto subcomplexes consisting of finite element subspaces formed with respect to arbitrary simplicial meshes. To keep the presentation as simple as possible we give details in the case of trimmed finite element spaces $\pol_r^-$, although the results will also work for the full polynomials spaces.  While there is overlap, there are also important differences between the work here and that in  \cite{falk2014local}. That paper makes use of weight functions $z_f^k$ \cite[pg. 2642]{falk2014local} which belong to the finite element spaces.  Here, instead, we use weight functions $\ZZ^k_r(\sig)$ that are not finite element functions. Their use allows us to avoid an extra correction step that seemed to be required in \cite{falk2014local}. In order to define $\ZZ^k_r(\sig)$ we use the formal adjoint of the exterior derivative and bubble functions to guarantee smoothness across interelement boundaries.  Our functions $\ZZ^k_r(\sig)$  rely on the existence of regular potentials for closed forms on contractible domains, for which we rely on the work
of Costabel and McIntosh \cite{costabel2010bogovskiui}. In particular, we use these results on the extended patch of a subsimplex. 

As in \cite{falk2014local}, we also first construct the projection onto the lowest-order space (e.g., the Whitney forms \cite{whitney2012geometric}). Then our projection for higher order elements uses the lowest order projection. An important difference is that we use alternative degrees of freedom for higher order elements. The degrees of freedom are essentially the ones used in the projection-based commuting interpolants developed by Demkowicz and collaborators \cite{demkowicz2005h1}. In fact, the degrees of freedom which we use are exactly the generalization of the degrees of freedom of the ones used be Melenk et al.  \cite{melenk2020commuting}. These degrees of freedom we allow us to define the higher order projections more efficiently. 

The paper is organized as follows. In the next section we give some preliminaries. In Section~\ref{loworder} we assume the existence of the weight functions $\ZZ^k_r(\sig)$ satifying certain requirements, and use them to construct the projections onto the lowest-order spaces, i.e., the spaces of the Whitney forms \cite{whitney2012geometric}. In the next section, we build on the lowest-order case to
construct the projection onto higher order finite elements, concluding with a statement of the main
result of the paper in Theorem~\ref{mainthm}.
Finally,  in Section \ref{Zsr} we give the deferred construction of the weight functions $\ZZ^k_r(\sig)$.

\section{Preliminaries}\label{prelim}

\subsection{Differential Forms}

The space of differential $k$-forms with smooth coefficients on a domain $S$ is denoted by $\La^k(S)$. The larger space allowing $L^2$ coefficients is denoted by $L^2 \La^k(S)$ and similarly $H^{\ell} \La^k(S)$ denotes the space of $k$-forms with coefficents in the Sobolev space $H^{\ell}(S)$.  The exterior derivative, denoted by $d^k$, maps $\La^k (S) \rightarrow \La^{k+1}(S)$ and extends to the spaces with less regularity. Finally, we define
\begin{equation*}
H\La^k(S):= \{\, u \in L^2 \Lambda^k(S)\,:\, d^k u \in L^2 \La^{k+1}(\Omega)\,\}. 
\end{equation*}

The Hodge star operator $\star$ maps $L^2\Lambda^k$ isomorphically onto $L^2\Lambda^{n-k}$ for each $k$.  Using it
we define the formal adjoint of $d^{k-1}$ by:
\begin{equation*}
 \de_k w= (-1)^k \star^{-1} d^{n-k}( \star w) \qquad \text{for } w \in \La^k(\Omega),
\end{equation*}
and the spaces
\begin{equation*}
H_{\delta}\La^k(S)=\{\, v \in L^2 \La^k(S)\,:\, \delta_k v \in L^2\La^{k-1}(S)\,\},
\quad \mathring{H}_{\delta}\La^k(S)=\{\, v \in H_{\delta}\La^k(S)\,:\, \text{tr}_{\partial S} \star v=0\,\},
\end{equation*}
the latter incorporating boundary conditions.
The adjoint relation between $d$ and $\delta$  may be expressed as
\begin{equation}\label{intparts}
\bl d^{k} \omega, v\br_S= \bl \omega, \de_{k+1} v\br_S \qquad \text{ for } \omega \in H\La^k(S), v\in  \mathring{H}_{\delta} \La^{k+1}(S),
\end{equation}
where $\bl\, \cdot\,,\, \cdot \,\br_S$ is the inner-product of $L^2 \Lambda^{k+1}(S)$.  We let $\| v\|_{L^2(S)}^2= \bl v,v\br_S$. We simply write $\bl\, \cdot\, ,\, \cdot \,\br$ when the domain $S$ is understood from the context.

\subsection{Simplicial complexes and co-boundaries} 
Let $\Omega \subset \mathbb{R}^n$ be a bounded domain and let $\Th$ be a simplicial triangulation of $\Omega$ consisting of $n$-simplices. We assume the shape regularity condition
\begin{equation*}
\frac{h_{\sig}}{\rho_{\sig}} \le C_S, \quad \sig \in \Th,
\end{equation*}
where $\rho_{\sig}$ is the diameter of the largest inscribed ball in $\sig$, $h_{\sig}$ is the diameter of $\sigma$, and $C_S>0$ is the shape regularity constant.  Closely related
to the triangulation $\Th$ is the associated simplicial complex $\del(\Th)$ consisting of all the simplices of $\Th$ and all their subsimplices of dimension $0$ through $n$.   We
denote by $\del_k(\Th)$, or simply $\del_k$ when the triangulation is clear, the collection of all the simplices in  $\del(\Th)$  of dimension $k$.  If $x_0,\ldots,x_k\in\mathbb{R}^n$ are the vertices of $\sig\in\del_k$, we may write
$[x_0,\ldots,x_k]$ for $\sig$, the closed convex hull of the vertices.
Often we need to endow a simplex with an orientation.  This is a choice of ordering of the vertices with two orders differing by an even permutation giving the same orientation.  If we select an ordering of all the vertices of $\Th$, this implies a default orientation for each of the simplices in $\del(\Th)$.

Let $x_1, \ldots, x_N$ be such an enumeration of the vertices of the mesh and let $\la_1, \ldots, \la_N$ be the continuous piecewise linear functions
such that $\la_i(x_j)=\delta_{ij}$. For  $\sigma =[x_{i_0} , \ldots, x_{i_n}] \in \del_n$
we define the bubble function
 \begin{equation*}
 b_{\sigma}:=\la_{i_0} \la_{i_1} \cdots \la_{i_n},
 \end{equation*}
a non-negative piecewise polynomial with support equal to $\sigma$.  To any simplex $\sigma =[x_{i_0} , \ldots, x_{i_k}] \in \del(\Th)$
we also associate the Whitney $k$-form $\phi_{\sig} \in H\La^k(\Omega)$ defined as
\begin{equation*}
\phi_{\sig}:=k! \sum_{j=0}^k (-1)^j \la_{i_j} d \la_{i_0} \w \cdots  \widehat{ d \la_{i_j}} \w \cdots  \w d \la_{i_k}.
\end{equation*}

For $\sigma \in \del(\Th)$ we define the \emph{star} of $\sigma$ as
\begin{equation*}
\st(\sigma)= \bigcup_{\tau \in \st_h(\sig)} \tau, \text{ where } \st_h(\sigma)= \{\, \tau \in  \del_n\,:\, \sigma \subset \tau \,\}, 
\end{equation*}
i.e., the union of all $n$-simplices containing $\sigma$. The \emph{extended star} of
$\sig$ is given by 
\begin{equation*}
\es(\sig)= \bigcup_{\tau \in \es_h(\sig)} \tau,  \text{ where } \es_h(\sig)=  \{\, \tau \in  \del_n\,:\, \sigma \cap  \tau \neq \emptyset \,\}, 
\end{equation*}
the union of $n$-simplices intersecting $\sigma$. As in \cite{ falk2014local}, we assume that $\es(\sig)$ is contractible for all $\sig$ in $\del(\Th)$, as is usually the case.

Associated with the simplicial complex are a chain complex
and cochain complex.  The space of $s$-chains is the vector space
\begin{equation*}
\C_k:= \{\, \sum_{ \sig \in \del_k} a_{\sig} \sigma\,:\,  a_\sig \in \mathbb{R} \,\}
\end{equation*}
where $\sig$ is given the default orientation and the same simplex with the opposite
orientation is identified with $-\sig$.
The boundary map $\partial_k:  \C_k \rightarrow \C_{k-1}$ is defined for $\sigma= [x_0, \ldots, x_k] \in \del_k$ by
\begin{equation*}
\partial_k \sigma=\sum_{j=0}^k (-1)^j [x_{0}, \ldots,\widehat{x_{j}}, \ldots, x_k].
\end{equation*}
The dual space of $\C_k$ is the space $\C^k$ of cochains. The basis $\del_k$ for
chains, leads to the dual basis $\sig^*$, $\sigma  \in \del_k$, for cochains, defined by
\begin{equation*}
\sig^*(\tau)=\delta_{\sig \tau}, \quad \sig, \tau\in\del_k
\end{equation*}
(using the Kronecker delta).
The coboundary operator $\dd^k: \C^k \rightarrow \C^{k+1}$ is defined by duality in the usual way:
\begin{equation}\label{defdd}
\dd^k X (\tau)=X(\partial_{k+1} \tau), \quad \tau \in \C_{k+1}, X \in \C^k.
\end{equation}
For the coboundary operator applied to a basis cochain we find that 
\begin{equation}\label{cob}
\dd^k [x_0 x_1 \ldots x_k]^*= \sum_{\substack{x\in \del_0 \\ [x \, x_0 x_1 \ldots x_k]  \in  \del_{k+1}}} [x \,  x_0 x_1 \ldots x_k]^*.
\end{equation}


\subsection{The FEEC forms}
 For integers $0\le k \le n$ and $r>0$, the space of trimmed polynomial $k$-forms of degree $r$ on $\R^n$  is
 $$
 \pol_{r}^-\La^k(\R^n)= \pol_{r-1}\La^k(\R^n)+\kappa \pol_{r-1}\La^{k+1}(\R^n)
 $$
 where $\kappa$ is the Kozul operator (see \cite{arnold2006finite}).  For a simplex $\tau$ in $\R^n$ of any dimension, the trimmed space on $\tau$ is given by restriction: $\pol_{r}^-\La^k(\tau)= \{\, \tr_\tau v\,:\, v \in  \pol_{r}^-\La^k(\R^n) \,\}$. Associated to any triangulation $\Th$ of $\R^n$ and to the integers $k$ and $r$ we then have the global trimmed finite element space $\sprk$, which is defined as
\begin{equation*}
\sprk= \{\, v \in H\La^k\,:\, v|_T \in \pol_{r}^-\La^k(T), \forall T \in \del_n \,\}.
\end{equation*}
The Whitney forms $\phi_{\sig}$, $\sig \in \del_k$, form a basis of $\Sp{1}{k}$.  The space $\sprk$ decomposes into a kernel portion and its orthogonal complement:
\begin{align*}
\Zz \sprk &:=  \{\, v \in \sprk \,:\, d^k v=0 \,\}, \\
\Zz^{\perp}  \sprk &:=   \{\, w \in  \sprk \,:\, \bl w, v\br=0, \forall v \in  \Zz \sprk \,\}.
\end{align*}

For a $k$-form $v$ that is smooth enough to admit an $L^1$ trace on some
$\sig\in\del(\Th)$, the de~Rham map defines the $k$-cochain $\dr^k v$ by
\begin{equation*}
\dr^k v (\sig)=\int_{\sig}\tr_\sig v,\quad \sig \in \del_k.
\end{equation*}
From Stokes theorem we easily see that
\begin{equation}\label{drcommute}
\dr^{k}(d^{k-1} v)(\sig)= \dr^{k-1} v(\partial_k \sig).
\end{equation}
The Whitney interpolant  $W^k: \C^k \rightarrow H\La^k(\Omega)$ is defined in term of
the Whitney forms by $W^k(\sig^*)=\phi_{\sig}$, so
\begin{equation*}
W^k(X)= \sum_{\sig \in \del_k} X(\sig) \phi_{\sig}.
\end{equation*}
The following properties of the Whitney interpolant are crucial  (see \cite[(4), (5) in page 139]{whitney2012geometric})
\begin{subequations}
\begin{alignat}{2}
\dr^k W^k X&= X,  \quad &&   X \in \C^k, \label{Wd0}\\
d^k W^k(X)&= W^{k+1}(\dd^k X), \quad  &&  X \in \C^k, \label{Wcommute} \\
\supp \phi_\sig &\subset  \st(\sig), &&  \sig \in \del_k. \label{Wsupport}
\end{alignat}
\end{subequations}
In particular, \eqref{Wd0} gives for $\sig \in \del_{s}$
\begin{equation}\label{Wd}
\dr^k \phi_\sig (\tau)=
\delta_{\sig\tau}, \quad \sig,\tau\in\del_k.
\end{equation}
It is easily shown that 
\begin{equation}\label{Wbound}
\|\phi_\sig\|_{L^2(\st(\sig))} \le C_W h_{\sig}^{\frac{n}{2}-k}, \quad  \sig \in \del_k,
\end{equation}
where $h_\sig$ is the local mesh size near $\sigma$.  (We may define $h_\sigma$ precisely as
the diameter of $\sigma$ if $s>0$ and as the diameter of $\st(\sig)$ if $\sig$ is a vertex.)

Now that we have introduced the Whitney and de Rham maps we can define the canonical projection, $\Pi_1^k$ onto the Whitney forms, which is given by
\begin{equation*}
\Pi_1^k:=W^k \dr^k.
\end{equation*}

Since we are assuming that $\es(\sig)$  is contractible for all $\sig \in \del(\Th)$, the local spaces  $\ssp{r}{\ell}{\es(\sigma)}$  form an exact sequence;  see \cite{arnold2018finite, arnold2006finite, falk2014local}.
\begin{proposition}
 Assume $\es(\sig)$ is contractible. For any $\sig \in \del(\Th)$ and $r \ge 1$ the following sequence is exact:
\begin{alignat}{4}\label{exact}
&\mathbb{R}
\stackrel{\subset}{\xrightarrow{\hspace*{0.5cm}}}\
 \ssp{r}{0}{\es_h(\sigma)}\
&&\stackrel{d^0}{\xrightarrow{\hspace*{0.5cm}}}\
 \ssp{r}{1}{\es_h(\sigma)}
&&\stackrel{d^1}{\xrightarrow{\hspace*{0.5cm}}}\
 \cdots
&&\stackrel{d^{n-1}}{\xrightarrow{\hspace*{0.5cm}}}\
 \ssp{r}{n}{\es_h(\sigma)}
\stackrel{}{\xrightarrow{\hspace*{0.5cm}}}\
0.
 \end{alignat}
 \end{proposition}

We will also need a discrete Poincar\'e inequality on the extended star $\es(\sigma)$.
\begin{proposition}
There exists a constant $C_P$ such that for all $\sig \in \del(\Th)$ one has
\begin{equation}\label{poincare}
 \|v\|_{L^2(\es(\sigma))} \le C_{P} h_{\sig} \|d^{\ell} v\|_{L^2(\es(\sigma))},  \quad v \in \Zz^{\perp}  \ssp{r}{\ell}{\es_h(\sigma)}.
\end{equation}
\end{proposition}
To prove this one uses the equivalence of norms on a finite dimensional space, together with a compactness argument and scaling by dilation. See \cite[Section 5.4]{arnold2006finite} and \cite[Section 5]{ falk2014local} for similar arguments.

The following proposition can be found Costabel and McIntosh \cite[Theorem 4.9 (c)]{costabel2010bogovskiui}.
It is proven using a generalized Bogovskii operator. 
\begin{proposition}\label{prop11}
Let $D$ be a bounded, contractible, Lipschitz domain.  Let $u \in \mathring{H}_{\delta}\La^k(D)$ satisfy $\delta_k u=0$, and also $\int_{D} u \, \text{vol}^n =0$ if $k=0$.  Then  there exists $\rho \in \mathring{H}^1 \La^{k+1}(D)$ such that $\delta_{k+1} \rho=u$. Moreover,
\begin{equation}\label{Ces}
|\rho|_{H^1(D)} \le C_D \|u\|_{L^2(D)}.  
\end{equation}
\end{proposition}
\noindent Using Friedrich's inequality (see \cite{gilbarg2015elliptic}) we have that $\|\rho\|_{L^2(D)}  \le C \text{diam}(D) C_D  \|u\|_{L^2(D)}$.  In  \cite{costabel2010bogovskiui} the authors do not track the constant $C_{D}$. However, in \cite{guzman2020estimation} it is shown that if $D$ is star-shaped with respect to a ball of similar diameter then the constant   $C_D$ can be bounded. Moreover, in \cite[Thm 33]{guzman2020estimation} bounds for the constants in slightly more general cases are given.  However, for arbitrary $\sig \in \del(\Th)$, the patch $\es(\sig)$ need not be star-shaped with respect to a ball and we cannot show that in general it satisfies the conditions of  \cite[Thm 33]{guzman2020estimation}. Therefore, we \emph{assume}  that the constants $C_{\es(\sigma)}$ are uniformly bounded.
\begin{assumption}\label{assump1}
Proposition \ref{prop11} is true when $D=\es(\sig)$ for all $\sig \in \del(\Th)$  with constants $C_{\es(\sigma)}$ uniformly bounded. 
\end{assumption}

Using that  $\text{diam}(\es(\sig)) \le C h_\sig$, we obtain the following result, which will use below. 
\begin{proposition}\label{bog}
Let $\sig \in \del(\Th)$. Assume the hypotheses of Proposition \ref{prop11} with $D=\es(\sig)$ and also  Assumption \ref{assump1}.  Then there exists a constant $C_{\delta}>0$  such that
\begin{equation}\label{deltabound}
\|\rho\|_{L^2(\es(\sig))} \le C_{\delta} h_{\sig} \|u\|_{L^2(\es(\sig))}.
\end{equation}
\end{proposition}

\section{Projection for the lowest order case $r=1$}\label{loworder}
In order to motivate our construction of an $L^2$-bounded projection, we recall the canonical  projection
which maps an element $u\in\La^k(\Omega)$ to
\begin{equation}\label{canonicalproj}
\Pi_1^k u=W^k \dr^k u=\sum_{\sig \in \del_k}  \dr^k u(\sig) \phi_\sig \in \Sp{1}{k}.
\end{equation}
In order that $\dr^k u(\sig)$ be well defined, $\tr_\sig u$ must be defined and integrable.   This is not the case for general $u\in L^2 \La^k(\Omega)$ when $k<n$. To obtain a projection that is well defined for  $u \in L^2 \La^k(\Omega)$, we replace $\dr^ku(\sig)$ with $\bl  \ZZ^k_r(\sig), u\br$ for a suitable $\ZZ^k_r(\sig)  \in L^2 \La^k(\Omega)$. (The subscript $r$, which refers to the polynomial degree, is introduced for the higher-order projections introduced in the next section.) In this section we
state the properties required of the differential form $\ZZ^k_r(\sig)$ and, assuming that such a form exists, develop an $L^2$-bounded projection into the Whitney forms.
We will verify the existence of a suitable form $\ZZ^k_r(\sig)$ in
Section~\ref{Zsr}.

Precisely, we shall show that for each  $r \ge 1$ and  $0\le k\le n$ there exist a linear operator $\ZZ^k_r: \C_k \rightarrow \mathring{H}_{\delta} \La^{k}(\Omega)$ which satisfies
\begin{subequations}\label{ZZall}
\begin{alignat}{2}
\bl \ZZ^k_r(\sig), u \br&= \,\dr^k u(\sig), \quad && u \in \sprk,\ \sig \in \del_k, \label{ZZ0r} \\
\delta_{k} \ZZ^k_r(\sig)&= \,\ZZ^{k-1}_r(\partial_k \sig) ,  \quad && \sig \in \del_k,  \label{ZZ2r}\\
\supp\ZZ^k_r(\sig)  &\subset \, \es(\sig),  &&    \sig \in \del_k, \label{ZZ3r} \\ 
\|\ZZ^k_r(\sig)\|_{L^2(\es(\sig))}  &\le \, C_Z h_{\sig} ^{-\frac{n}{2}+k},  \quad &&  \sig \in \del_k.   \label{ZZ4r}
\end{alignat}
\end{subequations}
\begin{definition}\label{pidef}
We define $\Prk{r}{k}: L^2 \La^k(\Omega) \rightarrow \Sp{1}{k}$ as 
\begin{equation*}
\Prk{r}{k} u:= \sum_{\sig \in \del_k} \bl  \ZZ^k_r(\sig), u\br \phi_\sig. 
\end{equation*}
\end{definition}
It follows directly from \eqref{ZZ0r}  and \eqref{canonicalproj} that $\Prk{r}{k}$ is an extension  of  $\Pi_1^k |_{\sprk}$:
\begin{lemma}
 For any $r \ge 1$, the  operator $\Prk{r}{k}: L^2 \La^k(\Omega) \rightarrow \Sp{1}{k}$ satisfies
 \begin{equation}
 \Prk{r}{k} u=\Pi_1^k u, \quad u \in  \sprk. \label{prksprk}
 \end{equation}
  
\end{lemma}
Moreover,  the operators $\Prk{r}{k}$ form bounded commuting projections:
\begin{theorem}\label{thmPi1r}
The operator  $\Prk{r}{k}: L^2 \La^k(\Omega) \rightarrow \Sp{1}{k}$ is a projection and the following commuting property holds:
\begin{equation}\label{commuting1r}
 d^k \Prk{r}{k} u= \Prk{r}{k+1} d^k u, \quad u \in H \La^k(\Omega).
\end{equation}
Moreover, if the mesh is shape-regular and Assumption \ref{assump1} holds, we obtain the following bound:
\begin{equation}\label{l2bound1r}
\|\Prk{r}{k} u\|_{L^2(T)} \le C \|u\|_{L^2(\es(T))},  \quad T \in \del_n,\ u\in L^2\La^k(\Omega).
\end{equation}
Finally, 
\begin{equation}\label{newPi1r}
\int_{\sig} \tr_{\sig} \Prk{r}{k} u=\int_{\sig} \tr_{\sig} u, \quad  \sig \in \del_k,
u \in \sprk.
\end{equation}
\end{theorem}

\begin{proof}
The fact that   $\Prk{r}{k}$ is a projection follows from \eqref{prksprk} and the fact that $\Pi_1^k$ is a projection.

To prove \eqref{l2bound1r}, let $T \in \del_n$. Since $\#\{\,\sig\in\del_k\,:\,\sig\subset T\,\}=c_1:=\binom{n+1}{k+1} $, we have 
\begin{align*}
\| \Prk{r}{k} u\|_{L^2(T)}^2&= \|\sum_{\substack{\sig \in \del_k \\ \sig \subset T}} \bl  \ZZ^k_r(\sig), u\br \phi_\sig\|_{L^2(T)}^2 
\le  c_1 \sum_{\substack{\sig \in \del_k \\
\sig \subset T}} |\bl  \ZZ^k_r(\sig), u\br|^2   \|\phi_\sig\|_{L^2(T)}^2 
\\
& \le c_1 \sum_{\substack{\sig \in \del_k \\
\sig \subset T}} \|u\|_{L^2(\es(\sig))}^2  \|\ZZ^k_r(\sig)\|_{L^2(\es(\sig))}^2   \|\phi_\sig\|_{L^2(T)}^2  \\
&\le c_1 C_W C_Z \sum_{\substack{\sig \in \del_k \\ \sig \subset T}} \|u\|_{L^2(\es(\sig))}^2 
\le c_1^2 C_W C_Z   \|u\|_{L^2(\es(T))}^2,
\end{align*}
where we used \eqref{Wbound},  \eqref{ZZ4r}. To prove \eqref{commuting1r} it suffices to prove 
\begin{equation}\label{815}
\dr^{k+1}(d^k \Prk{r}{k} u)(\tau) = \dr^{k+1}(\Prk{r}{k+1} d^k u)(\tau), \quad \tau \in \del_{k+1}.
\end{equation}

To this end, let $ \tau \in \del_{k+1}$ and use \eqref{Wd} to re-write the right-hand side as 
\begin{equation*}
\dr^{k+1}(\Prk{r}{k+1} d^k u)(\tau)= \sum_{\rho \in \del_{k+1}} \bl \ZZ^{k+1}_r(\rho), d^k u \br \dr^{k+1} \phi_\rho (\tau)
=\bl \ZZ^{k+1}_r(\tau), d^k u\br.
\end{equation*}
To treat the left-hand side we write
\begin{equation}\label{aux0}
 \partial \tau=\sig_0+ \cdots+\sig_{k+1},
\end{equation}
where the $\sig_i\in\del_k$ are the $k$-faces of $\tau$.
By \eqref{cob}, $\dd^k \sig_i^*$ is the sum of terms $\eta^*$ where $\eta$ runs over the 
$(k+1)$-simplices which contain $\sigma_i$ (taken with proper orientation), and, in particular, includes $\tau$.
Using again \eqref{Wd} we see that
\begin{equation}\label{aux2}
 \dr^{k+1}(W^{k+1} (\dd^k \sig_i^*))(\tau) =1, \quad  0 \le i \le k+1,
\end{equation}
while, if $\sig \in \del_k$ and $\sig$ is not contained in the boundary of $\tau$, then
\begin{equation}\label{aux3}
 \dr^{k+1}(W^{k+1} (\dd^k \sig^*))(\tau) =0.
\end{equation}

Therefore, 
 \begin{alignat*}{2}
 \dr^{k+1}(d^k \Prk{r}{k} u)(\tau) &= \sum_{\sig \in \del_k} \bl \ZZ^k_r(\sig), u \br \dr^{k+1} (d^k \phi_\sig)(\tau)  \qquad &&  \\
&= \sum_{\sig \in \del_k} \bl \ZZ^k_r(\sig), u \br \dr^{k+1} (W^{k+1} (\dd^k \sig^*))(\tau)  \qquad && \text{by } \eqref{Wcommute} \\
 &=  \sum_{i=0}^{k+1} \bl \ZZ^k_r(\sig_i), u \br  \qquad && \text{by }  \eqref{aux2}, \eqref{aux3} \\
&= \bl \ZZ^k_r(\partial_{k+1} \tau), u \br   \qquad && \text{by } \eqref{aux0} \\
&=  \bl \delta_{k+1} \ZZ^{k+1}_r(\tau), u \br && \text{by }  \eqref{ZZ2r} \\
&=  \bl \ZZ^{k+1}_r(\tau), d^k u \br && \text{by } \eqref{intparts}
 \end{alignat*}
 Thus, \eqref{815} holds.  Finally, \eqref{newPi1r} follows from \eqref{prksprk} and the definition of $\Pi_1^k$.
 \end{proof}
We see that $\pi_1^k:=\Prk{1}{k}$ is our desired projectionin the lowest-order case. In the next
section we will obtain the projection in the higher order case $r>1$ as a correction to
$\Prk{r}{k}$.

\section{Higher-order Elements}\label{highorder}
\subsection{Idea of the construction}
Next we discuss the strategy for constructing the projection in the general case.  The first step is to decompose the space $\sprk$ using the projection $\Pi_1^k$.  For each $r\ge 1$ we have
\begin{equation*}
\sprk=\Pi_1^k \sprk \oplus (I-\Pi_1^k) \sprk=\Sp{1}{k} \oplus M_r^k,
\end{equation*}
where we have set $M_r^k= (I-\Pi_1^k) \sprk$.  Note that
\begin{equation}\label{defMrk}
M_r^k=\{ v \in \sprk: \Pi_1^k v=0 \}= \{\,v \in \sprk \,:\, \int_{\sig} \tr_{\sig} v=0,\forall \sig \in \del_k\,\}.  
\end{equation}
In particular, $M_1^k=0$.  Also, using Stokes theorem we easily see that the spaces $M_r^k$ with $r$ fixed and $k$ increasing form a sub-complex of the complex formed by the $\sprk$.  The key step is to construct a projection $Q_r^k: L^2 \La^k(\Omega) \rightarrow M_r^k$ that is local, $L^2$-bounded and commutes with the exterior derivative
\begin{equation*}
 d^k Q_r^k u= Q_r^{k+1} d^k u, \quad  u \in H \La^k(\Omega).
\end{equation*}
Then we define $\pi_r^k: L^2 \La^k(\Omega) \rightarrow \sprk$ for all $r \ge 1$ as
\begin{equation}\label{defproj}
\pi_r^k u: =\Prk{r}{k} u+Q_r^k(u-\Prk{r}{k} u), \quad  u \in L^2 \La^k(\Omega). 
\end{equation}
If $u \in \sprk$ then $u-\Prk{r}{k} u \in M_r^k$ by \eqref{newPi1r} and hence $\pi_r^k u=u$,  so $\pi_r^k$ is indeed a projection. Moreover, one can easily show that it commutes with the exterior derivative.

\subsection{Alternative degrees of freedom for $v \in \pol_r^{-} \La^k(\tau)$}
We now turn to the key step of constructing the projection $Q_r^k$.  For this, it is useful to use degrees of freedom (dofs) for the space $\sprk$ different than the canonical degrees of freedom described in \cite{arnold2006finite}. Instead we will use dofs developed by Demkowicz and collaborators \cite{demkowicz2005h1}, a generalization of the ones found in Melenk et al.~\cite{melenk2020commuting}.

Let $\tau$ be any simplex and consider the polynomial differential form spaces   $\pol_r^{-} \La^k(\tau)$   and $\mathring{\pol}_r^{-} \La^k(\tau)= \{\, v \in  \pol_r^{-} \La^k(\tau)\,:\, \tr_{\partial \tau} v=0 \,\}$ where $k \le \dim \tau$.   We have the following  exact sequence
\begin{alignat}{6}\label{exactPr0}
&0
\stackrel{}{\xrightarrow{\hspace*{0.5cm}}}\
\mathring{\pol}_r^-\La^0(\tau)\
&&\stackrel{d}{\xrightarrow{\hspace*{0.5cm}}}\
\mathring{\pol}_r^-\La^1(\tau)
&&\stackrel{d}{\xrightarrow{\hspace*{0.5cm}}}\
 \cdots
%
%
&&\stackrel{d}{\xrightarrow{\hspace*{0.5cm}}}\
\mathring{\pol}_r^-\La^{\dim \tau}(\tau)  
&&\stackrel{\int_{\tau}}{\xrightarrow{\hspace*{0.5cm}}}\
\R
&&\stackrel{}{\xrightarrow{\hspace*{0.5cm}}}\
0.
 \end{alignat}
Letting
\begin{equation*}
\polb_r^{-} \La^k(\tau)=
\begin{cases}\mathring{\pol}_r^{-} \La^k(\tau), &k < \dim \tau,\\
\{\, v \in  \mathring{\pol}_r^{-} \La^k(\tau)\,:\, \int_{\tau} v=0\,\}, &  k=\dim\tau,
\end{cases}
\end{equation*}
we obtain the exact sequence:
\begin{alignat}{5}\label{exactPrcup}
&0
\stackrel{}{\xrightarrow{\hspace*{0.5cm}}}\
\polb_r^-\La^0(\tau)\
&&\stackrel{d}{\xrightarrow{\hspace*{0.5cm}}}\
\polb_r^-\La^1(\tau)
&&\stackrel{d}{\xrightarrow{\hspace*{0.5cm}}}\
 \cdots
%
%
&&\stackrel{d}{\xrightarrow{\hspace*{0.5cm}}}\
\polb_r^-\La^{\dim \tau}(\tau)
&& \stackrel{}{\xrightarrow{\hspace*{0.5cm}}}\
0
 \end{alignat}

Next, we decompose  $\mathring{\pol}_r^{-} \La^k(\tau)$  into the kernel of $d$ and the space orthogonal to the kernel:
\begin{equation}\label{717}
\mathring{\pol}_r^{-} \La^k(\tau)=  \Zz \mathring{\pol}_r^{-} \La^k(\tau)\oplus \Zz^{\perp} \mathring{\pol}_r^{-} \La^k(\tau) ,
\end{equation}
where
\begin{alignat*}{1}
\Zp^k(\tau) &:= \{\, z \in \mathring{\pol}_r^{-} \La^k(\tau)\,:\,  d z=0  \,\}, \\
\Zpp^k(\tau)&:= \{\, \eta \in \mathring{\pol}_r^{-} \La^k(\tau)\,:\, \bl \eta, z \br_{\tau}=0, \forall z \in \Zp^k(\tau) \,\}.
\end{alignat*}
Note that $\Zp^k(\tau)= \mathring{\pol}_r^{-} \La^k(\tau)$ when $\dim \tau=k$. Similarly, $\polb_r^{-} \La^k(\tau)$ decomposes into $\ZZp^k(\tau)$ and $\ZZpp^k(\tau)$. 

We know that  \cite[Theorem~4.14]{arnold2006finite}
\begin{equation}\label{716}
 \dim \mathring{\pol}_r^{-} \La^k(\tau)= \dim  \pol_{r+k-\dim \tau-1} \La^{\dim \tau -k}(\tau).
 \end{equation}
 In particular, $\mathring{\pol}_r^{-} \La^k(\tau)=0$ if $\dim \tau < k$ or if $\dim \tau \ge r+k$. Therefore, by  \cite[Theorem~4.13]{arnold2006finite} that 
 \begin{equation}\label{917}
  \dim \pol_r^{-} \La^k(T)=\sum_{\tau \in \del(T)}  \dim \mathring{\pol}_r^{-} \La^k(\tau).
 \end{equation}

 In order to introduce the dofs efficiently we define the bilinear form 
 \begin{equation*}
 \bll u, v \brr_{\tau}= \bl \Q_{\tau} u, \Q_{\tau}  v \br_{\tau}+ \bl du, dv \br_{\tau}, \quad  u,v  \in H \Lambda^k(\tau), 
 \end{equation*}
 where $\Q_{\tau}$ is the $L^2$-orthogonal projection onto  $\Zp^k(\tau)$ given by
 \begin{equation*}
 \bl \Q_{\tau} u,   w \br_{\tau}=\bl u, w\br_{\tau}, \qquad w \in \Zp^k(\tau).
 \end{equation*}
 Note that   $\bll \,\cdot\,,\, \cdot \,\brr_{\tau}$ is an inner product on the space  $\mathring{\pol}_r^{-} \La^k(\tau)$.

We now give the dofs for  $\pol_r^{-} \La^k(T)$ and prove their unisolvence. 
\begin{lemma}\label{lemmadofZ}
Let $T$ be a simplex. Then, $w \in  \pol_r^{-} \La^k(T)$ is determined by
\begin{alignat}{1}
\bll \tr_{\tau} w, y \brr_{\tau}, &  \quad y \in  \mathring{\pol}_r^{-} \La^k(\tau), \, \tau \in \del(T).  \label{dofZ}
\end{alignat}
\end{lemma}
Note that \eqref{dofZ} is vacuous unless $k\le \dim \tau < r+k$. 
\begin{proof}
By \eqref{917} the total number of dofs in \eqref{dofZ} is the same as the dimension of $\pol_r^{-} \La^k(T)$. Suppose that the dofs \eqref{dofZ} of $w$ vanish. We must show that $w=0$. We can do this by induction on $\dim T$. The base case $\dim T=0$ is trivial.  By the induction step $\tr_{\tau} w=0$ for all   $\tau \subset  \del(T)$ with $\tau  \neq  T$ which in particular implies that  $w \in \mathring{\pol}_r^{-} \La^k(T)$. Thus, choosing $\tau=T$ and $y=w$ in \eqref{dofZ} gives that $\bll w, w \brr_{T}=0$. Since  $\bll \,\cdot\,,\, \cdot \,\brr_{T}$ is an inner-product on $\mathring{\pol}_r^{-} \La^k(T)$, this implies that $w = 0$ 
\end{proof}
As a corollary we immediately obtain dofs for the global finite element space $\sprk$.
\begin{corollary}\label{corollaryVunique}
A differential $k$ form $w \in  \sprk$ is uniquely determined by
\begin{alignat*}{1}
\bll \tr_{\tau} w, y \brr_{\tau}, &  \quad y \in  \mathring{\pol}_r^{-} \La^k(\tau),\, \tau \in \del(\Th).  
\end{alignat*}
\end{corollary}
\begin{remark}\label{remark}
Let $\dim \tau=k$. Then the volume form of $\tau$, $\text{vol}_{\tau}$, belongs to $\mathring{\pol}_r^{-} \La^k(\tau)$ for all $r \ge 1$ and thus $\int_{\tau} w= \bl w, \text{vol}_{\tau} \br_{\tau} = \bll \tr_{\tau} w,  \text{vol}_{\tau} \brr_{\tau}$ is always  a dof of  $w \in  \sprk$. If $r=1$ then these are the only dofs given in Corollary  \ref{corollaryVunique} and they coincide with the canonical dofs in the case $r=1$. 
\end{remark}
In fact, based on this remark we have the following corollary.
\begin{corollary}\label{corollaryMunique}
A differential $k$ form $w \in  M_r^k$ is uniquely determined by
\begin{alignat*}{1}
\bll \tr_{\tau} w, y \brr_{\tau}, &  \quad y \in  \polb_r^{-} \La^k(\tau),\, \tau \in \del(\Th).  
\end{alignat*}
\end{corollary}

\subsection{Discrete Extensions}
In this subsection we define some key spaces and extension operators. A differential form $w\in \sprk$ is determined by the dofs given in Corollary~\ref{corollaryVunique}.  For $\sig\in\del(\Th)$, we define $G_r^k(\sig)$ as the space of all $w\in \sprk$ for which all those dofs vanish except those associated to the simplex
$\sigma$.  We note that $G_r^k(\sig)=0$ if $\dim\sig < k$ or $\dim\sig \ge r+k$.  In any case,
if $w\in G_r^k(\sig)$, then $\supp w\subset \st(\sig)$.

A simple consequence of Lemma \ref{lemmadofZ} is the following:
\begin{lemma}
Let $v \in G_r^{k}(\sig)$, and suppose that $\tau \in \del(\Th)$ does not contain $\sigma$.  Then,
\begin{equation}\label{lemmatrfree}
\tr_{\tau} v=0.
\end{equation}
In particular, this occurs when $\dim \tau < \dim \sig$ or $\dim\tau=\dim\sig$ but $\tau \neq \sig$. 
\end{lemma}

Next we define an operator from $E_{\sig}: H \La^k(\sig) \rightarrow G_r^{k}(\sig)$ as follows. For any $\rho \in H \La^k(\sig)$ let $E_{\sig} \rho \in G_r^{k}(\sig)$ satisfy
\begin{equation} \label{E}
\bll \tr_{\sig} E_{\sig} \rho , y \brr_{\sig}=\bll \rho , y \brr_{\sig},   \quad   y \in   \mathring{\pol}_r^{-} \La^k(\sig).
\end{equation}
Note that $E_{\sig}$ maps a $k$-form on $\sig$ to a piecewise polynomial $k$-form on $\Omega$.
In view of  Lemma \ref{lemmadofZ}, we see that 
\begin{equation}\label{trEsig}
v=E_{\sig} \tr_{\sigma}v, \quad v \in  G_r^{k}(\sig).
\end{equation}
The next result shows that the operator $E_{\sig}$ is an extension operator if we restrict ourselves to  $\mathring{\pol}_r^{-} \La^k(\sig)$ and that it commutes the exterior derivative if we further restrict ourselves to  $\polb_r^{-} \La^k(\sig)$. 
\begin{lemma}
Let $\sig \in \Delta(\Th)$, $0 \le k \le n$ and  $r \ge 1$. Then,
\begin{alignat}{2}
\tr_{\sig} E_{\sig} \rho&=\rho,  \quad  &&  \rho \in  \mathring{\pol}_r^{-} \La^k(\sig), \label{trE} \\
d E_{\sig} \rho&= E_{\sig} d \rho,  \quad  && \rho \in  \polb_r^{-} \La^k(\sig). \label{commuteE}
\end{alignat}
\end{lemma}

\begin{proof}
We prove  \eqref{trE} first. Using \eqref{lemmatrfree} we have that $\tr_{\sig} E_{\sig} \rho  \in \mathring{\pol}_r^{-} \La^k(\sig)$. Then, $\phi=\tr_{\sig} E_{\sig} \rho-\rho \in \mathring{\pol}_r^{-} \La^k(\sig)$ and it satisfies
\begin{alignat*}{1}
\bll \phi , y \brr_{\sig}=0, &  \quad   y \in \mathring{\pol}_r^{-} \La^k(\sig).
\end{alignat*}
From this we conclude that $\phi \equiv 0$ which proves \eqref{trE}.

Next we turn to the proof of \eqref{commuteE}.  Let  $\rho \in  \polb_r^{-} \La^k(\sig)$ and set $w=d E_{\sig} \rho- E_{\sig} d \rho  \in \Sp{r}{k+1}$. We will show that 
\begin{alignat}{1}
\bll \tr_{\tau} w, y \brr_{\tau} = 0, &  \quad  y \in \mathring{\pol}_r^{-} \La^{k+1}(\tau), \, \tau \in \del(\Th). \label{8716}
\end{alignat}
This in turn implies that $w \equiv 0$ by Lemma \ref{dofZ}.
We  consider separately the cases $\tau \neq \sig$ and $\tau=\sig$.

{\bf{Case 1 : $\tau  \neq \sig$}}.   Since $E_{\sig} d \rho \in G_r^{k+1}(\sig)$ and $\tau \neq \sig$ we have
\begin{alignat*}{1}
\bll \tr_{\tau} w, y \brr_{\tau} = \bll \tr_{\tau} d E_{\sig} \rho , y \brr_{\tau}= \bl \tr_{\tau} d  E_{\sig} \rho, \Q_{\tau} y \br_{\tau},
\end{alignat*}
where we used that $d \circ d=0$.  We only need to consider the case $\dim \tau  \ge k+1$ since otherwise $\mathring{\pol}_r^{-} \La^{k+1}(\tau)=0$.  If $\dim \tau >k+1$ then by the exact sequence \eqref{exactPr0} we have that $ \Q_{\tau} y=dm$ for some $m \in \Zpp^k(\tau)$. Hence, 
\begin{equation}\label{aux781}
 \bl \tr_{\tau} d  E_{\sig} \rho, \Q_{\tau} y \br_{\tau}= \bl \tr_{\tau} d  E_{\sig} \rho,d m \br_{\tau}= \bll  \tr_{\tau} E_{\sig} \rho, m \brr_{\tau}=0,
\end{equation}
where we used that $\Q_{\tau} m=0$ and then that $E_{\sig} \rho \in G_r^k(\sig)$. On the other hand,  if $\dim \tau=k+1$,  then we chose $c \in \R$ so that $\int_{\tau} (\Q_{\tau} y-c \, \vol_{\tau})=0$. Thus, again by the exactness of exact sequence \eqref{exactPr0} we have $\Q_{\tau} y- c \, \vol_{\tau} =d m$ for some $m \in \Zpp^k(\tau)$. Hence, using \eqref{aux781} we have 
\begin{equation*}
 \bl \tr_{\tau} d  E_{\sig} \rho, \Q_{\tau} y \br_{\tau}= \bl \tr_{\tau} d  E_{\sig} \rho,d m \br_{\tau}+ c \int_{\tau} \tr_{\tau} d  E_{\sig} \rho=   c \int_{\tau} \tr_{\tau} d  E_{\sig} \rho.
\end{equation*}
Then, if we write $\partial \tau= \sum_{i=1}^\ell \eta_i$ where $\eta_i \in \del_k(\Th)$  we have using Stokes formula 
\begin{equation*}
 \int_{\tau} \tr_{\tau} d  E_{\sig} \rho=  \sum_{i=1}^\ell  \int_{\eta_i} \tr_{ \eta_i} E_{\sig} \rho.
\end{equation*}
If $\eta_i \neq \sig$ then $\int_{\eta_i} \tr_{ \eta_i} E_{\sig} \rho=0$ since  $E_{\sig} \rho \in G_r^k(\sig)$ where we used Remark \ref{remark}. If $\eta_i=\sig$  then using, \eqref{E}, and Remark \ref{remark} we get  $\int_{\sig} \tr_{\sig} E_{\sig} \rho=\int_{\sig} \tr_{\sig} \rho=0$ where we used our hypothesis that  $\rho \in  \polb_r^{-} \La^k(\sig)$. Thus, we have shown \eqref{8716} in the case $\tau \neq \sig$.

{\bf{Case 2 : $\tau  = \sig$}}. Using the definition of $E_{\sig}$, \eqref{E}, and \eqref{trE} we have
\begin{alignat*}{1}
\bll \tr_{\sig} w, y \brr_{\sig} =   \bll \tr_{\sig}  (d E_{\sig} \rho- E_{\sig} d \rho) , y \brr_{\sig} =\bll d (\tr_{\sig }E_{\sig} \rho-  \rho) , y \brr_{\sig}=0.
\end{alignat*}
\end{proof}

In order to define the projections it is helpful to identify an orthonormal basis for  $\polb_r^{-} \La^k(\sig)$ which we will denote by $\mathfrak{p}_r^k(\sig)$. For $0 \le k \le \dim \sig$ we let  $\mathfrak{z}_r^{k,\perp}(\sig)$ be a basis of  $\Zz^{\perp} \mathring{\pol}_r^{-} \La^k(\sig)$ satisfying
\begin{equation}\label{perp0}
\bll p, q\brr_{\sig}= \delta_{pq}, \qquad  p, q  \in \mathfrak{z}_r^{k,\perp}(\sig).
\end{equation}
 For $1 \le k \le \dim \sig$,  we define $\mathfrak{z}_r^{k}(\sig):= \{\, dp\,:\, p \in \mathfrak{z}_r^{k-1,\perp}(\sig) \,\}$ and  if $k=0$ define  $\mathfrak{z}_r^{0}(\sig)=\emptyset$. Finally, we define 
 \begin{equation}\label{mathfrakp}
 \mathfrak{p}_r^k(\sig):=\mathfrak{z}_r^{k}(\sig) \cup \mathfrak{z}_r^{k,\perp}(\sig).
 \end{equation}
  That this is a basis of $\polb_r^{-} \La^k(\sig)$  follows from the exactness of \eqref{exactPr0}.  Indeed, it is an orthonormal basis with respect to the inner product $\bll \,\cdot\,,\, \cdot\, \brr_{\sig}$:
 \begin{equation}\label{perp0_2}
\bll p, q\brr_{\sig}= \delta_{pq}, \qquad  p, q  \in  \mathfrak{p}_r^k(\sig).
\end{equation}

We can use these orthonormal functions to give a representation formula for any  $M_r^k$. This follows from \eqref{perp0_2},  Corollary \ref{corollaryMunique}, the definition of $E_{\sig}$ and \eqref{trE}.
\begin{lemma}
Any $u \in M_r^k$ can be uniquely written as 
\begin{equation}\label{lemma112}
u=  \sum_{\sig \in \del(\Th)} \sum_{g \in \mathfrak{p}_r^k(\sig)} \bll \tr_{\sig} u, g \brr_{\sig} E_{\sig} g.
\end{equation}
\end{lemma}

\subsection{The projection $Q_r^k$}
As usual, let $0\le k\le n$ and $r\ge 1$ be integers.  The projection $Q_r^k: L^2 \La^k(\Omega) \rightarrow M_r^k$ will be defined in terms of linear operators
\begin{equation*}
\UU_r^k(\tau): \polb_r^{-} \La^k(\tau)\to \mathring{H}_{\delta}\La^k(\Omega),\quad
g\mapsto \UU_r^k(\tau,g),
\end{equation*}
to be defined for each $\tau \in \Delta(\Th)$.  To define $\UU_{r}^k(\tau)$ we make use of the space  $M_r^k(\st_h(\tau))=\{\, v|_{\st(\tau)}\,:\, v \in M_r^k \,\}$, defined by \eqref{defMrk} with the triangulation $\Th$ replaced by the subtriangulation $\st_h(\tau)$.  We will also use
the function
\begin{equation}\label{defBBb}
\BBb_{\tau}:= \sum_{\sig \in \st_h(\tau)} b_{\sig} ,
\end{equation}
which is the superposition of the bubble functions on the $n$-simplices comprising $\st_h(\tau)$. Clearly it
is supported in $\st(\tau)$ and vanishes on any simplex in $\Th$ of dimension less than $n$.

Turning to the definition of $\UU_r^k(\tau)$, we
note that if $\dim \tau <k$, then $\polb_r^{-} \La^k(\tau)$ vanishes and so $\UU_r^k(\tau)=0$. For $k \le \dim \tau$ we define
$\UU_r^k(\tau)$ separately on $\ZZp^k(\tau)$ and $\ZZpp^k(\tau)$, namely
we define
\begin{equation*}
\UU_{r}^k(\tau, g)=\BBb_{\tau} \beta,\quad g \in \ZZp^k(\tau),   
\end{equation*}
where $\beta \in M_r^k(\st_h(\tau))$ is the unique solution to
\begin{equation}\label{aux763}
\bl \BBb_{\tau} \beta, u \br_{\st(\tau)}= \bl g, \tr_{\tau} u\br_{\tau}, \quad  u \in M_r^k(\st_h(\tau)), 
\end{equation}
and
\begin{equation*}
\UU_{r}^k(\tau, g)= \delta_{k+1} \UU_{r}^{k+1}(\tau, d g),\quad g \in \ZZpp^k(\tau).   
\end{equation*}

The following lemma establishes the properties of $\UU_r^k(\tau)$.
\begin{lemma}
For every $\tau \in \del(\Th)$ the following properties hold: 
\begin{subequations}\label{UUall}
\begin{alignat}{2}
 \bl  u,  \UU_r^k(\tau, g)\br&=\,  \bll   \tr_{\tau} u, g \brr_{\tau},  \quad  &&   g  \in \polb_r^{-} \La^k(\tau), u \in M_r^k, \label{UU2-1}  \\
 \delta_{k+1} \UU_{r}^{k+1} (\tau, dg )&=\, \UU_r^{k} (\tau, g),  \, \quad &&   g \in \ZZpp^k(\tau) \label{UU3-2}, \\
 \supp  \UU_r^k(\tau, g)  &\subset \, \st(\tau),  && g  \in \polb_r^{-} \La^k(\tau),   \label{UU3}  \\
  \|\UU_{r}^k(\tau, dm)\|_{L^2(\st(\tau))}  &\le \,  C_U h_{\tau}^{(\dim \tau-n)/2}, \, \quad &&   m \in \mathfrak{z}_r^{k-1,\perp}(\tau).  \label{UU4}
\end{alignat}
\end{subequations}
\end{lemma}

\begin{proof}
 We see that \eqref{UU3-2} and \eqref{UU3} follow immediately from the definition of $\UU$.  Let us prove \eqref{UU2-1} first in the case $g \in \ZZp^k(\tau)$. In this case $dg=0$ and $g= \Q_{\tau} g$ and so
 \begin{alignat*}{2}
 \bl  u,  \UU_r^k(\tau, g)\br=   \bl g, \tr_{\tau} u\br_{\tau}= \bl  \Q_{\tau} g,  \Q_{\tau} (\tr_{\tau} u)\br_{\tau} = \bll \tr_{\tau} u, g \brr_{\tau},
 \end{alignat*}
 where we used \eqref{aux763}. On the other hand, suppose that $g \in \ZZpp^k(\tau)$. Then we note that $dg \in  \ZZp^{k+1}(\tau)$ and so by the previous case 
$\bl  d u,  \UU_r^k(\tau, dg)\br= \bll \tr_{\tau} du, dg \brr_{\tau}$ for any  $u \in M_r^k$ since $du \in M_r^{k+1}$. Thus, for $u \in M_r^k$ we can use the definition of $\UU$ and integration by parts to get 
\begin{alignat*}{1}
 \bl  u,  \UU_r^k(\tau, g)\br=   \bl  u,  \delta_{k+1} \UU_r^k(\tau, dg)\br=  \bl  d u,  \UU_r^k(\tau, dg)\br  =  \bll \tr_{\tau} du, dg \brr_{\tau}. 
 \end{alignat*}
 Using that $d \circ d g=0$, $\Q_{\tau} dg= dg$ and $\Q_{\tau} g=0$  we get
\begin{equation*}
 \bll \tr_{\tau} du, dg \brr_{\tau} = \bl \Q_{\tau}  (\tr_{\tau} du),   \Q_{\tau} dg \br_{\tau} =   \bl d \tr_{\tau} u,   dg \br_{\tau} =    \bll \tr_{\tau} u,   g \brr_{\tau} .
\end{equation*}
This proves \eqref{UU2-1}. The estimate \eqref{UU4} follows from a scaling argument. 
\end{proof}

 Now we define the projection $Q_r^k: L^2 \La^k(\Omega) \rightarrow M_r^k$ as
\begin{equation*}
Q_{r}^k u:= \sum_{\sig \in \del(\Th)} \sum_{g \in \mathfrak{p}_r^k(\sig)} \bl u, \UU_r^k(\sig,g) \br E_{\sig} g, \quad  u \in L^2 \La^k(\Omega).
\end{equation*}

Before proving the main result of this section we will need the following estimate that follows from the definition of $E_{\tau}$ and a scaling argument
\begin{equation}\label{scaleE}
 \frac{1}{h_{\tau}}\|E_{\tau} (m)\|_{L^2(\st(\tau))} + \|E_{\tau} (d m)\|_{L^2(\st(\tau))} \le h_{\tau}^{(n-\dim \tau)/2},  \quad m \in \mathfrak{z}_r^{k-1,\perp}(\tau). 
\end{equation}


\begin{lemma}\label{lemmapir}
The operator  $Q_r^k: L^2 \La^k(\Omega) \rightarrow M_r^k$ is a projection and the following commuting property holds
\begin{equation}\label{commuting1}
 d^k Q_r^k u= Q_r^{k+1} d^k u, \quad u \in H \La^k(\Omega).
\end{equation}
Moreover, 
\begin{equation}\label{l2boundpi}
\|Q_{r}^k u\|_{L^2(T)} \le C \|u\|_{L^2(\es(T))},  \quad T \in \del_n.
\end{equation}
\end{lemma}

\begin{proof}
The fact that  $Q_r^k$  is a projection follows from \eqref{lemma112}, \eqref{UU2-1}. To prove \eqref{commuting1} we use \eqref{commuteE}  and \eqref{mathfrakp} to get
\begin{equation*}
d Q_{r}^k u=  \sum_{\sig \in \del(\Th)}  \sum_{g \in \mathfrak{p}_r^{k}(\sig)} \bl  \UU_{r}^k(\sig, g) , u \br E_{\sig}(d g)= \sum_{\sig \in \del(\Th)}  \sum_{g \in \mathfrak{z}_r^{k, \perp}(\sig)} \bl    \UU_{r}^k(\sig, g), u \br E_{\sig}(d g).
\end{equation*}
On the other hand, using integration integration by parts we obtain
\begin{equation*}
Q_{r}^{k+1} d u=  \sum_{\sig \in \del(\Th)}  \sum_{m \in \mathfrak{p}_r^{k+1}(\sig)} \bl  \UU_{r}^{k+1}(\sig, m) , d u\br E_{\sig} m =  \sum_{\sig \in \del(\Th)}  \sum_{m \in \mathfrak{p}_r^{k+1}(\sig)} \bl \delta  \UU_{r}^{k+1}(\sig, m) , u\br E_{\sig} m.
\end{equation*}
If use \eqref{UU3-2} then we see that $\delta \UU_{r}^{k+1}(\sig, m)=0$ when $m \in  \mathfrak{z}_r^{k+1, \perp}(\sig)$. Hence, 
\begin{equation*}
Q_{r}^{k+1} d u= \sum_{\sig \in \del(\Th)}  \sum_{m \in \mathfrak{z}_r^{k+1}(\sig)} \bl \delta  \UU_{r}^{k+1}(\sig, m) , u\br E_{\sig} m =  \sum_{\sig \in \del(\Th)}  \sum_{g \in \mathfrak{z}_r^{k, \perp}(\sig)} \bl \delta \UU_r^{k+1}(\sig,dg) , u\br E_{\sig} (dg).
\end{equation*}
We now see that \eqref{commuting1} follows from another application of \eqref{UU3-2}. Finally, one can easily establish \eqref{l2boundpi} using \eqref{scaleE} and \eqref{UU4}. 
\end{proof}

\subsection{The final projection}
Having defined $Q_r^k$, the projection $\pi_r^k:L^2 \La^k(\Omega)\rightarrow\sprk$ is defined by
\eqref{defproj}.  As pointed out there, it is indeed a projection operator.  Moreover, since
$M_1^k=0$, the operator $Q_1^k$ vanishes and so $\pi_1^k$ coincides with $\Prk{1}{k}$.
We end this section by noting that $\pi_r^k$ is indeed an $L^2$ bounded projection operator which commutes with the exterior derivative and which is local.  This is the main result of the paper.
\begin{theorem}\label{mainthm}
The operator $\pi_r^k:  L^2 \La^k(\Omega) \rightarrow \sprk$ is a projection. It commutes with the exterior derivative: 
\begin{equation}\label{commutePi}
 d^k \pi_r^k u= \pi_r^{k+1} d^k u, \quad u \in H \La^k(\Omega).
\end{equation}
Moreover, if the mesh is shape regular and Assumption \ref{assump1} holds, then we have the
local $L^2$ estimates
\begin{subequations}\label{l2boundr}
\begin{alignat}{2}
\|\pi_{1}^k u\|_{L^2(T)} &\le C \|u\|_{L^2(\es(T))},  \quad  && T \in \del_n, \\
\|\pi_{r}^k u\|_{L^2(T)} &\le C \|u\|_{L^2(\es^2(T))},  \quad && T \in \del_n, r \ge 2,
\end{alignat}
\end{subequations}
where 
\begin{equation*}
\es^2(\sig)= \bigcup_{\substack{T \subset \es(\sig)\\ T \in \del_n}} \es(T).
\end{equation*}

\end{theorem}

\begin{proof}
From \eqref{commuting1} and \eqref{commuting1r} we obtain the commutativity:
\begin{equation*}
 d^k \pi_r^k u= d^k \Prk{r}{k} u+ d^k Q_r^k(u-\Prk{r}{k} u) =\Prk{r}{k+1} d^k u+ Q_r^{k+1}(d^k u- \Prk{r}{k+1} d^k u) =\pi_r^{k+1} d^k u.
\end{equation*}
The result \eqref{l2boundr} follows from \eqref{l2boundpi} and \eqref{l2bound1r}.
\end{proof}

\section{Construction of $\ZZ^k_r$} \label{Zsr}

It remains to construct the linear operators $\ZZ^k_r: \C_k \rightarrow \mathring{H
}_{\delta} \La^{k}(\Omega)$ satisfying the properties \eqref{ZZ0r}-\eqref{ZZ4r} whose
existence was asserted in Section~\ref{loworder}. We will
again use a superposition of bubble functions as in \eqref{defBBb}, but now defined with respect
to the extended star of a simplex:
\begin{equation*}
\bb_{\tau}= \sum_{\sig \in \es_h(\sigma)} b_{\sig} .
\end{equation*}
Thus $\bb_{\sig}$ is supported in $\es(\sigma)$ and vanishes on any simplex in $\Th$ of dimension less than $n$. 
\begin{lemma}\label{lemmaBs}
Let $\sig \in  \del_k$ and let $L:\Zz^{\perp} \ssp{r}{k}{\es_h(\sig)} \rightarrow \mathbb{R}$ be a linear functional. Then the following problem has a  unique solution:

Find $v \in \Zz^{\perp}  \ssp{r}{k}{\es_h(\sig)}$ satisfying
\begin{equation}
\bl \bb_{\sig} d^k v, d^k u\br=L(u),  \quad u \in   \Zz^{\perp}  \ssp{r}{k}{\es_h(\sig)}.
\end{equation}
\end{lemma}
\begin{proof}
This is a square linear system and therefore we only have to prove uniqueness. Suppose that $\bl \bb_{\sig} d^k u, d^k u\br=0$ for some $u \in \Zz^{\perp}  \ssp{r}{k}{\es_h(\sig)}$. Then, given the property of $\bb_{\sig}$ we have that for each $\tau \in \del_n$ with $\tau \subset \es(\sig)$ we have that  $d^k u$ vanishes on $\tau$. Hence, $d^k u $ vanishes on $\es(\sig)$ or $u \in \Zz  \ssp{r}{k}{\es_h(\sig)}$. Thus, $u$ must be zero. 
\end{proof}

We also need the following bound that follows from inverse estimates:
\begin{equation}\label{dr0}
\dr^k(u)(\sig) \le C h_{\sig}^{-\frac{n}{2}+k}\|u\|_{L^2(\es(\sig))}, \quad  u \in   \ssp{r}{k}{\es_h(\sig)}.
\end{equation}

Finally, we will need the following result which follows from a scaling argument:
\begin{equation}\label{bubble1}
\|\sqrt{\bb_{\sig}} d \rho\|_{L^2(\es(\sig))} \le C \|\bb_{\sig} d\rho\|_{L^2(\es(\sig))}, \quad  \rho \in  \ssp{r}{k}{\es_h(\sig)}.  
\end{equation}

\begin{theorem}
Assume that the mesh is shape regular and assume that Assumption \eqref{assump1} holds. For $0 \le k \le n$,  there exists a linear operator $\ZZ^k_r: \C_k \rightarrow \mathring{H}_{\delta} \La^{k}(\Omega)$ satisfying \eqref{ZZ0r}--\eqref{ZZ4r}. 
\end{theorem}

\begin{proof}
The proof is by induction on $k$. To initialize we need to define
$\ZZ^0_r(\sig)\in  \mathring{H}_{\delta} \La^{k}(\Omega)$ for $\sig\in \del_0$ and $r\ge 1$.
The $0$-simplex $\sig$ consists of a single vertex, say $p$.  Let
$\eta=\chi_{\es(\sigma)}/|\es(\sigma)|$ denote the characteristic function of $\es(\sigma)$ normalized to have integral unity. Hence, $\bl u, \eta \br=u(p)$ for any constant function $u$.
Invoking Lemma~\ref{lemmaBs}, we define $v \in \Zz^{\perp} \ssp{r}{0}{\es_h(\sig)}$ by
\begin{equation}\label{605}
\bl \bb_{\sig} d^0 v, d^0 u\br=u(p) -\bl \eta, u \br, \quad  u \in   \Zz^{\perp}  \ssp{r}{0}{\es_h(\sig)}.
\end{equation}
Note that $\Zz \ssp{r}{0}{\es_h(\sig)}$ consists of constant functions, so both the left-hand side and right-hand side of the equation in \eqref{605} vanish for $u \in  \Zz \ssp{r}{0}{\es_h(\sig)}$. Hence
\begin{equation}\label{805}
\bl \bb_{\sig} d^0 v, d^0 u\br=u(p)-\bl \eta, u \br, \quad  u \in   \ssp{r}{0}{\es_h(\sig)}.
\end{equation}
We now define
\begin{equation}\label{defZ0}
\ZZ^0_r(\sig):=\eta+ \delta_1\big(\bb_{\sig} d^0 v\big).
\end{equation}

To finish the initial step of the induction we verify that the operator $\ZZ^0_r$ has the desired properties \eqref{ZZ0r}--\eqref{ZZ4r}. 
For $u \in \Sp{r}{0}$, we have by \eqref{defZ0}, \eqref{intparts}, and \eqref{805} that
\begin{equation*}
\bl \ZZ^0_r(\sig), u \br= \bl \eta+\delta_1\big( \bb_{\sig} d^0 v\big), u \br = \bl \eta, u \br+ \bl  \bb_{\sig} d^0 v, d^0 u \br = u(p),
\end{equation*}
which is \eqref{ZZ0r} in the case $k=0$.  For $k=0$ both sides of \eqref{ZZ2r} trivially vanish,
and the locality condition \eqref{ZZ3r} is clear.
Finally, to prove \eqref{ZZ4r} we note that 
\begin{equation}\label{871}
\|\eta\|_{L^2(\es(\sig))} \le C h_{\sig}^{-n/2}.
\end{equation}
Moreover, by \eqref{805} we have 
\begin{alignat*}{2}
\|\sqrt{\bb_{\sig}} d^0 v\|_{L^2(\es(\sig))}^2 &\le |v(p)|+  \|\eta\|_{L^2(\es(\sig))} \|v\|_{L^2(\es(\sig))} \qquad && \\
&\le h_{\sig}^{-\frac{n}{2}} \|v\|_{L^2(\es(\sig))} \qquad && \text{ inverse estimate and } \eqref{871}  \\
&\le h_{\sig}^{-\frac{n}{2}+1} \|d^0 v\|_{L^2(\es(\sig))} \qquad && \text{ by } \eqref{poincare}.
\end{alignat*}
Hence, using \eqref{bubble1} we get 
\begin{equation*}
\|\sqrt{\bb_{\sig}} d^0 v\|_{L^2(\es(\sig))} \le C h_{\sig}^{-\frac{n}{2}+1}. 
\end{equation*}
By another inverse estimate we have
\begin{equation*}
\|\delta_1\big(\bb_{\sig} d^0 v\big)\|_{L^2(\es(\sig))} \le C h_{\sig}^{-\frac{n}{2}}. 
\end{equation*}
This combined with \eqref{871}  shows \eqref{ZZ4r} for the case $k=0$.

To complete the induction, we suppose that we have constructed $\ZZ^\ell_r$ satisfying \eqref{ZZ0r}--\eqref{ZZ4r} for any $\ell <k$, and we construct $\ZZ^k_r$.  Let $\sig \in \del_k$ be arbitary.  By Proposition \ref{prop11} we have the existence of $\eta \in \mathring{H}_{\delta} \La^k(\es(\sig))$ such that 
\begin{equation}\label{743}
\delta_k \eta=\ZZ^{k-1}_r(\partial_k \sigma) \quad \text{ on } \es(\sig).
\end{equation}
Here we used that $\delta_{k-1}\ZZ^{k-1}_r(\partial_k \sigma)=\ZZ^{k-2}_r(\partial_{k-1} \partial_k \sigma)=0$ if $k \ge 2$ which follows from our induction hypothesis and \eqref{ZZ2r} and also $\bl \ZZ^0_r (\partial_1 \sigma), 1 \br_{\es(\sig)}=0$ if $k=1$.

We then define
\begin{equation}\label{defZZs}
\ZZ^k_r(\sig):=\eta+ \delta_{k+1}(\bb_{\sig} d^k v) ,
\end{equation}
where $v \in \Zz^{\perp} \ssp{r}{k}{\es_h(\sig)}$ solves (by Lemma \ref{lemmaBs})
\begin{equation}\label{331}
\bl \bb_{\sig} d^k v, d^k u\br=\dr^k u(\sig) -\bl \eta, u \br, \quad  u \in  \Zz^{\perp} \ssp{r}{k}{\es_h(\sig)}.
\end{equation}
Note if $u \in   \Zz \ssp{r}{k}{\es_h(\sig)}$ then the left-hand side of \eqref{331} vanishes. Moreover, by Proposition \ref{exact} there exists $w \in \ssp{r}{k-1}{\es_h(\sig)}$ such that $d^{k-1} w=u$. Hence, 
\begin{alignat*}{2}
\dr^k u(\sig)  -\bl \eta, u \br&= \dr^k d^{k-1} w (\sig)  -\bl \eta, d^{k-1} w \br \\
&= \dr^{k-1} w(\partial \sig)   -\bl \delta_k \eta,  w \br \qquad && \text{ by } \eqref{drcommute}, \eqref{intparts} \\
&=  \dr^{k-1} w(\partial \sig)  -\bl \ZZ^{k-1}_r(\partial \sigma),  w \br \qquad && \text{ by } \eqref{743} \\
&=0 \qquad && \text{ by  our induction hypothesis}. 
\end{alignat*}
Thus, we have 
\begin{equation}\label{331-2}
\bl \bb_{\sig} d^k v, d^k u\br=\dr^k u(\sig) -\bl \eta, u \br, \quad  u \in  \ssp{r}{k}{\es_h(\sig)}.
\end{equation}

We complete the induction by proving that indeed $\ZZ^k_r(\sig)$ has the desired properties. By \eqref{743}, \eqref{defZZs} and the fact that $\delta_k \delta_{k+1}=0$ we have that \eqref{ZZ2r} holds. Clearly \eqref{ZZ3r}  holds.  We next prove \eqref{ZZ0r} holds.  Let $u \in \sprk$,  then
\begin{alignat*}{2}
\bl \ZZ^k_r(\sig), u \br&= \bl \eta+\delta_{k+1}( \bb_{\sig} d^k v, u\br \qquad && \text{ by } \eqref{defZZs} \\
&= \bl \eta, u \br+ \bl  \bb_{\sig} d^k v, d^k u \br \qquad &&    \text{ by } \eqref{intparts} \\
&= \dr^ku(\sig) \qquad && \text{ by } \eqref{331-2}.
\end{alignat*}
To prove \eqref{ZZ4r} we note by \eqref{deltabound} that
\begin{equation*}
\|\eta\|_{L^2(\es(\sig))} \le C_{\delta} h_{\sig} \|\ZZ^{k-1}_r(\partial \sig)\|_{L^2(\es(\sig))}.
\end{equation*}
By our induction hypothesis \eqref{ZZ4r} we have 
\begin{equation*}
 \|\ZZ^{k-1}_r(\partial \sig)\|_{L^2(\es(\sig))} \le C \,h_{\sig}^{-\frac{n}{2}+k-1}.
 \end{equation*}
 Thus, 
 \begin{equation}\label{r512}
\|\eta\|_{L^2(\es(\sig))} \le C h_{\sig}^{-\frac{n}{2}+k}.
 \end{equation}

Using \eqref{331}   we get 
\begin{alignat*}{2}
\|\sqrt{\bb_{\sig}} d^k v\|_{L^2(\es(\sig))}^2 &\le  |\dr^k v|+\|\eta\|_{L^2(\es(\sig))} \|v\|_{L^2(\es(\sig))} \qquad &&  \\
&\le C  h_{\sig}^{-\frac{n}{2}+k} \|v\|_{L^2(\es(\sig))}  \qquad && \text{ by }  \eqref{r512}, \eqref{dr0}\\
&\le C  h_{\sig}^{-\frac{n}{2}+k+1 }  \|d^k v\|_{L^2(\es(\sig))}.  \qquad && \text{ by }  \eqref{poincare} 
\end{alignat*}
Therefore, using \eqref{bubble1} we have 
\begin{equation*}
\|\sqrt{\bb_{\sig}} d^k v\|_{L^2(\es(\sig))} \le C  h_{\sig}^{-\frac{n}{2}+k+1 }.
\end{equation*}
If we now use an inverse estimate we get 
\begin{alignat}{1}\label{r513}
\|\delta_{k+1}\big(\bb_{\sig} d^k v \big)\|_{L^2(\es(\sig))} \le  \frac{C}{h_\sig}\|\bb_{\sig} d^k v \big\|_{L^2(\es(\sig))} \le  \frac{C}{h_\sig}\|\sqrt{\bb_{\sig}} d^k v \big\|_{L^2(\es(\sig))} \le C h_{\sig}^{-\frac{n}{2}+k}.
\end{alignat}
Combining \eqref{r512} and \eqref{r513} gives \eqref{ZZ4r}.

\end{proof}

\bibliographystyle{abbrv}
\bibliography{references}

\end{document}